\documentclass[11pt,a4paper,twoside]{article}
\usepackage{amssymb,amsfonts}
\usepackage{fancyhdr}
\usepackage{amsthm}
\usepackage{hyperref}
\usepackage{color}
\usepackage[T1]{fontenc}

\usepackage{afterpage}  %

\usepackage{amssymb,amsfonts,amsthm}

\newcounter{thrm}
\newcounter{lmm}
\newcounter{prp}
\newcounter{crll}
\newcounter{rmrk}
\newenvironment{Theorem}[1][\hspace{-1ex}]%
 {\par\addvspace{0.3em}\noindent\refstepcounter{thrm}\bf Theorem~\thethrm\hspace{1ex}{\rm #1}.\,~\it}%
 {\rm\par\addvspace{0.3em}}

\newenvironment{Lemma}[1][\hspace{-1ex}]%
 {\par\addvspace{0.3em}\noindent\refstepcounter{lmm}\bf Lemma~\thelmm\hspace{1ex}{\rm #1}.\,~\it}%
 {\rm\par\addvspace{0.3em}}

\newenvironment{Prop}[1][\hspace{-1ex}]%
 {\par\addvspace{0.3em}\noindent\refstepcounter{prp}\bf Proposition~\theprp\hspace{1ex}{\rm #1}.\,~\it}%
 {\rm\par\addvspace{0.3em}}

\newenvironment{Corol}[1][\hspace{-1ex}]%
 {\par\addvspace{0.3em}\noindent\refstepcounter{crll}\bf Corollary~\thecrll\hspace{1ex}{\rm #1}.\,~\it}%
 {\rm\par\addvspace{0.3em}}

\newenvironment{Note}[1][\hspace{-1ex}]%
 {\par\addvspace{0.3em}\noindent\refstepcounter{rmrk}\bf Remark~\thermrk\hspace{1ex}{\rm #1}.\,~\rm}%
 {\rm\par\addvspace{0.3em}}

\newenvironment{Proof}[1][\hspace{-1ex}]%
 {\par\addvspace{0.3em}\noindent{\em Proof\hspace{1ex}#1.\ }}%
 {\hfill$\Box$\par\addvspace{0.3em}}

\def\eqdf{\stackrel{\mathrm{\scriptscriptstyle def}}=}
\def\ij#1#2{_{#1}^{#2}}

\begin{document}

\afterpage{\rhead[]{\thepage}
 \chead[\small D. S. Krotov, V. N. Potapov, P. V. Sokolova]%
 {\small  On reconstructing reducible $n$-ary quasigroups and switching subquasigroups}
 \lhead[\thepage]{}
}

\title{On reconstructing reducible $n$-ary quasigroups\\
       and switching subquasigroups}
\author{Denis S. Krotov, Vladimir N. Potapov, Polina V. Sokolova}
\renewcommand\today{}
\maketitle
\begin{abstract}
{
(1) We prove that, provided $n\geq 4$, a permutably reducible $n$-ary quasigroup
is uniquely specified by its values on the $n$-ples containing zero.
(2) We observe that for each $n,k\geq 2$ and $r\leq \lfloor k/2\rfloor$ there exists
a reducible $n$-ary quasigroup of order $k$ with an $n$-ary subquasigroup of order $r$.
As corollaries, we have the following:
(3) For each $k\geq 4$ and $n\geq 3$ we can construct
a permutably irreducible $n$-ary quasigroup of order $k$.
(4) The number of $n$-ary quasigroups of order $k>3$ has double-exponential growth as $n\to\infty$;
it is greater than $\exp\exp(n\ln \lfloor k/3\rfloor)$ if $k\geq 6$, and
$\exp\exp({\ln 3 \over 3}n - 0.44)$ if $k=5$.
}
\end{abstract}
\footnote{The paper will appear in the \href{http://www.quasigroups.eu}{Quasigroups And Related Systems}, 16 (2008) no.1}

\section[Introduction]{\hspace{-1em}.~\,Introduction}

An $n$-ary operation $f:\Sigma^n\to \Sigma$,
where $\Sigma$ is a nonempty set,
is called
an \emph{$n$-ary qua\-si\-group} or \emph{$n$-qua\-si\-group}
(\emph{of order $|\Sigma|$}) iff
in the equality $z_{0}=f(z_1, \ldots , z_n)$ knowledge of any $n$ elements
of $z_0$, $z_1$, \ldots , $z_n$ uniquely specifies the remaining one \cite{Belousov}.

An $n$-ary quasigroup $f$ is \emph{permutably reducible} iff
$$f(x_1,\ldots,x_n)=h\left(g(x_{\sigma (1)},\ldots,x_{\sigma (k)}),
\linebreak[1]
x_{\sigma (k+1)},\ldots,x_{\sigma (n)}\right)$$
where $h$ and $g$ are $(n-k+1)$-ary and $k$-ary quasigroups, $\sigma$ is a permutation,
and $1<k<n$.
In what follows we omit the word ``permutably''
because we consider only such type of reducibility.

We will use the following standard notation: $x\ij{i}{j}$ denotes $x_i,x_{i+1},\ldots, x_j$.

In Section~\ref{s:reduc}
we show that
a reducible $n$-qua\-si\-group can be reconstructed by its values
on so-called `shell'.
`Shell' means the set of variable values with at least one zero.

In Section~\ref{s:subq}
we consider the questions of imbedding
$n$-qua\-si\-groups of order $r$ into $n$-qua\-si\-groups of order $k\geq 2r$.

In Section~\ref{s:irreduc}
we prove that for all $n\geq 3$ and $k\geq 4$
there exists an irreducible $n$-qua\-si\-group of order $k$.
Before, the question of existence of irreducible $n$-qua\-si\-groups
was considered by Belousov and Sandik \cite{BelSan} ($n=3$, $k=4$),
Frenkin \cite{Frenkin} ($n\geq 3$, $k=4$), Borisenko \cite{Borisenko} ($n\geq 3$, composite finite $k$),
Akivis and Goldberg \cite{Gold75,Gold76,AkGo} (local differentiable $n$-qua\-si\-groups),
Glukhov \cite{Glu76} ($n\geq 3$, infinite $k$).

In Sections~\ref{s:number}
and~\ref{s:comp}
we prove the double-exponential (of type $\exp\exp(c(k)n)$) lower bound on
the number $|Q(n,k)|$ of $n$-qua\-si\-groups of finite order $k\geq 4$.
Before, the following asymptotic results on the number of
$n$-qua\-si\-groups of fixed finite order $k$ were known:
\begin{itemize}
\item $|Q(n,2)|=2$.
\item $|Q(n,3)|=3\cdot 2^n$, see, e.g., \cite{LaywineMullen};
a simple way to realize this fact is to show by induction that
the values on the shell uniquely specify an $n$-qua\-si\-group of order $3$.
\item $|Q(n,4)|= 3^{n+1}2^{2^n +1}(1+o(1))$ \cite{PotKro:asymp,KroPot:Lyap}.
\end{itemize}
Note that by the ``number of $n$-qua\-si\-groups'' we mean the number of mutually
different $n$-ary qua\-si\-group operations $\Sigma^n\to\Sigma$
for a fixed $\Sigma$, $|\Sigma|=k$ (sometimes, by this phrase one means
the number of isomorphism classes).
As we will see, for every $k\geq 4$ there is $c(k)>0$ such that $|Q(n,k)|\geq 2^{2^{c(k)n}}$.
More accurately (Theorem~\ref{th:Num}),
if $k=5$ then $|Q(n,5)|\geq 2^{3^{n/3-const}}$;
for even $k$ we have $|Q(n,k)|\geq 2^{(k/2)^n}$;
for $k\equiv 0\bmod 3$ we have $|Q(n,k)|\geq 2^{n(k/3)^n}$;
and for every $k$ we have $|Q(n,k)|\geq 2^{1.5\lfloor k/3\rfloor^n}$.
Observe that dividing by the number (e.g., $(n+1)!(k!)^n$) of any natural equivalences
(isomorphism, isotopism, paratopism,\ldots)
does not affect these values notably;
so, for the number of equivalence classes almost the same bounds are valid.
For the known exact numbers of $n$-qua\-si\-groups of order $k$ with small values of $n$ and
$k$, as well as the numbers of equivalence classes for different equivalences, see the
recent paper of McKay and Wanless \cite{MK-W:small}.

\section[On reconstructing reducible $n$-quasigroups]{{\hspace{-1em}.~\,On reconstructing reducible $n$-quasigroups}}\label{s:reduc}
In what follows the constant tuples $\bar o$, $\bar \theta$ may be
considered as all-zero tuples.
From this point of view,
the main result of this section states that a reducible $n$-qua\-si\-group
is uniquely specified by its values on the `shell',
where the `shell' is
the set of $n$-ples with at least one zero.
Lemma~\ref{p:12} and its corollary concern the case
when the groups of variables in the decomposition
of a reducible $n$-qua\-si\-group are fixed.
In Theorem~\ref{th:vosst} the groups of variables are not specified;
we have to require $n\geq 4$ in this case.

\begin{Lemma}[(a representation of a reducible $n$-qua\-si\-group
by the superposition of retracts)]\label{p:12}\mbox{}
Let $h$ and $g$ be an
$(n-m+1)$-
and $m$-qua\-si\-groups, let $\bar o\in\Sigma^{m-1}$, $\bar \theta\in\Sigma^{n-m}$, and let
\begin{eqnarray}\nonumber
&f(x,\bar y, \bar z)\eqdf h(g(x,\bar y),\bar z),&
\\
h_0(x,\bar z)\eqdf f(x,\bar o,\bar z),
&
g_0(x,\bar y)\eqdf f(x,\bar y,\bar \theta),
&
\delta(x)\eqdf f(x,\bar o,\bar \theta)
                                \label{e:7}
\end{eqnarray}
where
$x\in\Sigma$,
$\bar y\in\Sigma^{m-1}$,
$\bar z\in\Sigma^{n-m}$.
Then
\begin{equation}
f(x,\bar y, \bar z)\equiv h_0(\delta^{-1}(g_0(x,\bar y)),\bar z).
                                \label{e:8}
\end{equation}
\end{Lemma}
\begin{Proof}
It follows from (\ref{e:7}) that
$$
h_0(\cdot,\bar z)\equiv h(g(\cdot,\bar o),\bar z),
\quad
g_0(x,\bar y)\equiv h(g(x,\bar y),\bar \theta),
\quad
\delta^{-1}(\cdot)\equiv g^{-1}(h^{-1}(\cdot,\bar \theta),\bar o).
$$
Substituting these representations of
$h_0$, $g_0$, $\delta^{-1}$
to (\ref{e:8}), we can readily
verify its validity.
\nolinebreak
\end{Proof}

\begin{Corol}\label{l:sync}
Let $q_{in},q_{out},f_{in},f_{out}:\Sigma^2\to\Sigma$ be qua\-si\-groups,
$q \eqdf q_{out}(x_1,q_{in}(x_2,x_3))$, $f \eqdf f_{out}(x_1,f_{in}(x_2,x_3))$, and $(o_1,o_2,o_3) \in \Sigma^3$.
Assume that for all $(x_1,x_2,x_3) \in \Sigma^3$ it holds
$$
q(o_1,x_2,x_3) = f(o_1,x_2,x_3),\quad
q(x_1,o_2,x_3) = f(x_1,o_2,x_3).
$$
Then $q(\bar x)= f(\bar x) $ for all $\bar x \in \Sigma^3$.
\end{Corol}
\begin{Theorem}\label{th:vosst}
Let $q,f:\Sigma^n\to\Sigma$ be reducible $n$-qua\-si\-groups, where $n\geq 4$;
and let $o\ij{1}{n} \in \Sigma^n$.
Assume that for all $i\in \{1,\ldots,n\}$ and for all $ x\ij{1}{n} \in \Sigma^n$ it holds
\begin{equation}\label{e:obol}
q(x\ij{1}{i-1}, o_i, x\ij{i+1}{n}) = f(x\ij{1}{i-1}, o_i, x\ij{i+1}{n}).
\end{equation}
Then $q( x\ij{1}{n})= f( x\ij{1}{n}) $ for all $ x\ij{1}{n} \in \Sigma^n$.
\end{Theorem}
\begin{Proof}
(*) We first proof the claim for $n=4$.
Without loss of generality (up to coordinate permutation and/or interchanging $q$ and $f$), we can assume that
one of the following holds for some quasigroups $q_{in},q_{out},f_{in},f_{out}$:

Case 1) $q(x\ij{1}{4})=q_{out}(x_1,q_{in}(x_2,x_3,x_4))$, $f(x\ij{1}{4})=f_{out}(x_1,f_{in}(x_2,x_3,x_4))$;

Case 2) $q(x\ij{1}{4})=q_{out}(x_1,q_{in}(x_2,x_3,x_4))$, $f(x\ij{1}{4})=f_{out}(x_1,f_{in}(x_2,x_3),x_4)$;

Case 3) $q(x\ij{1}{4})=q_{out}(x_1,q_{in}(x_2,x_3),x_4)$, $f(x\ij{1}{4})=f_{out}(x_1,f_{in}(x_2,x_3),x_4)$;

Case 4) $q(x\ij{1}{4})=q_{out}(x_1,q_{in}(x_2,x_3,x_4))$, $f(x\ij{1}{4})=f_{out}(f_{in}(x_1,x_2,x_3),x_4)$;

Case 5) $q(x\ij{1}{4})=q_{out}(x_1,q_{in}(x_2,x_3,x_4))$, $f(x\ij{1}{4})=f_{out}(f_{in}(x_1,x_4),x_2,x_3)$;

Case 6) $q(x\ij{1}{4})=q_{out}(x_1,x_2,q_{in}(x_3,x_4))$, $f(x\ij{1}{4})=f_{out}(x_1,f_{in}(x_2,x_3),x_4)$;

Case 7) $q(x\ij{1}{4})=q_{out}(x_1,q_{in}(x_2,x_3),x_4)$, $f(x\ij{1}{4})=f_{out}(f_{in}(x_1,x_4),x_2,x_3)$.

1,2,3) Take an arbitrary $x_4$ and denote
$q'(x_1,x_2,x_3) \linebreak[2] \eqdf q(x_1,x_2,x_3,x_4)$
and
$f'(x_1,x_2,x_3) \linebreak[2] \eqdf f(x_1,x_2,x_3,x_4)$.
Then, by Corollary~\ref{l:sync},
we have $q'(\bar x) = f'(\bar x)$ for all $\bar x\in \Sigma^3$; this proves the statement.

4) Fixing $x_4:=o_4$ and applying (\ref{e:obol}) with $i=4$, we have
$f_{out}(f_{in}(x_1,x_2,x_3),o_4)=q_{out}(x_1,q_{in}(x_2,x_3,o_4))$,
which leads to the
representation
$f_{in}(x_1,x_2,x_3) = h_{out}(x_1,h_{in}(x_2,x_3))$
where
$h_{out}(x_1,\cdot) \eqdf f_{out}^{-1}(q_{out}(x_1,\cdot),o_4)$ and
$h_{in}(x_2,x_3) \eqdf q_{in}(x_2,x_3,o_4)$.
Using this representation, we find that $f$ satisfies the condition of Case 2) for
some $f_{in},f_{out}$.
So, the situation is reduced to the already-considered case.

5) Fixing $x_4:=o_4$ and using (\ref{e:obol}), we obtain the
decomposition $f_{out}(\cdot,\cdot,\cdot) = h_{out}(\cdot,h_{in}(\cdot,\cdot))$
for some $h_{in},h_{out}$. We find that
$q$ and $f$ satisfy the conditions of Case 2).

6) Fixing $x_4:=o_4$ and using (\ref{e:obol}), we get the
decomposition $q_{out}(\cdot,\cdot,\cdot) = h_{out}(\cdot,h_{in}(\cdot,\cdot))$.
Then, we again reduce to Case 2).

7) Fixing $x_4:=o_4$ we derive
the decomposition $f_{out}(\cdot,\cdot,\cdot) = h_{out}(\cdot,h_{in}(\cdot,\cdot))$,
which leads to Case 3).

(**) Assume $n>4$.
It is straightforward to show that we always can choose four indexes $1\leq i<j<k<l\leq n$
such that for all
$x\ij{1}{i-1}$,
$x\ij{i+1}{j-1}$,
$x\ij{j+1}{k-1}$,
$x\ij{k+1}{l-1}$,
$x\ij{l+1}{n}$
the $4$-qua\-si\-groups
\begin{eqnarray}\nonumber
q'_{x\ij{1}{i-1}x\ij{i+1}{j-1}x\ij{j+1}{k-1}x\ij{k+1}{l-1}x\ij{l+1}{n}}(x_i,x_j,x_k,x_l) \eqdf q(x\ij{1}{n}),
\\ \nonumber
f'_{x\ij{1}{i-1}x\ij{i+1}{j-1}x\ij{j+1}{k-1}x\ij{k+1}{l-1}x\ij{l+1}{n}}(x_i,x_j,x_k,x_l) \eqdf f(x\ij{1}{n})
\end{eqnarray}
are reducible. Since these $4$-qua\-si\-groups satisfy the hypothesis of the lemma,
they are identical, according to (*). Since they coincide for every values of the parameters,
we see that $q$ and $f$ are also identical.
\end{Proof}

\begin{Note}
If $n=3$ then the claim of Lemma~\ref{th:vosst} can fail.
For example, the reducible $3$-qua\-si\-groups
$q(x\ij{1}{3})\eqdf (x_1 * x_2) * x_3$
and
$f(x\ij{1}{3})\eqdf x_1 * (x_2 * x_3)$ where $*$ is a binary quasigroup
with an identity element $0$ $($i.\,e., a loop$)$
coincide if $x_1=0$, $x_2=0$, or $x_3=0;$ but they are not identical
if $*$ is nonassociative.
\end{Note}

\section[Subquasigroup]{{\hspace{-1em}.~\,Subquasigroup}}\label{s:subq}
Let $q:\Sigma^n\to \Sigma$ be an $n$-quasirgoup and $\Omega \subset\Sigma$.
If $g=q|_{\Omega ^n}$ is an $n$-quasirgoup then we will say that
$g$ is a \emph{subquasigroup} of $q$ and $q$ is \emph{$\Omega$-closed}.


\begin{Lemma}\label{l:subq}
For each finite $\Sigma$ with $|\Sigma|=k$ and $\Omega \subset\Sigma$
with $|\Omega |\leq \lfloor k/2\rfloor$
there exists a reducible 
$n$-qua\-si\-group $q:\Sigma^n\to \Sigma$
with a subquasigroup $g:\Omega^n\to \Omega$.
\end{Lemma}
\begin{Proof}
By Ryser theorem on completion of a Latin $s\times r$ rectangular up to
a Latin $k\times k$ square ($2$-qua\-si\-group) \cite{Ryser51},
there exists a $\Omega $-closed $2$-qua\-si\-group $q:\Sigma^2\to \Sigma$.

To be constructive, we suggest a direct formula for the case
$\Sigma=\{0,\ldots,k-1\}$,
$\Omega =\{0,\ldots,r-1\}$ where $k\geq 2r$ and $k-r$ is odd:
$$
\begin{array}{rcll}
  q_{k,r}(i,j) &=& (i+j)\bmod r,\quad &i<r, j<r; \nonumber\\
  q_{k,r}(r+i,j) &=& (i+j)\bmod (k-r)+r,\quad &j<r; \nonumber\\
  q_{k,r}(i,r+j) &=& (2i+j)\bmod (k-r)+r,\quad &i<r; \nonumber\\
  q_{k,r}(r+i,r+j)&=&\multicolumn{2}{l}
  { \cases{
  (i-j)  \bmod  (k-r)\quad &if $(i-j)\bmod (k-r)<r$,\cr
  (2i-j)  \bmod  (k-r)+r\quad &otherwise.}
  }
\end{array}
$$
In the following four examples the second and the fourth value arrays correspond to
$q_{5,2}$ and $q_{7,2}$:
\begin{equation}\label{e:57}
\def\arraystretch{0.8}
\arraycolsep=0.6ex
\def\rsrs#1{\raisebox{-1.5pt}{{#1}}}
\def\0{\rsrs{0}}
\def\1{\rsrs{1}}
\def\2{\rsrs{2}}
\def\3{\rsrs{3}}
\def\4{\rsrs{4}}
\def\5{\rsrs{5}}
\def\6{\rsrs{6}}
{4{:}}~
\begin{array}{|cccc|}
\hline
{\0}  & \multicolumn{1}{c|}{\1}  & \2  & \3\\
{\1}  & \multicolumn{1}{c|}{\0}  & \3  & \2\\
\cline{1-2}
\2  & \3  & {\0}  & \1\\
\3  & \2  & {\1}  & \0\\
\hline
\end{array}
\qquad
{5{:}}~
\begin{array}{|cc|ccc|}
\hline
{\0}  & {\1}  & \2  & \3&\4\\
{\1}  & {\0}  & \3  & \4&\2\\ \hline
\2  & \4  & {\0}  & \1&\multicolumn{1}{|c|}{\6}\\
\cline{3-3}\cline{5-5}\3  & \2  & \4  &\multicolumn{1}{|c}{\0}&{\1}\\
\cline{3-3}\cline{4-4}\4  & \3  & \1  & \multicolumn{1}{|c}{\5}&\multicolumn{1}{|c|}{\0}\\
\hline
\end{array}
\qquad
{6{:}}~
\begin{array}{|cccccc|}
\hline
{\0}  & \multicolumn{1}{c|}{\1}  & \2  & \3 & \4  & \5\\
{\1}  & \multicolumn{1}{c|}{\0}  & \3  & \2 & \5  & \4\\
\cline{1-2}
\4 & \5 & \0 & \1 & \2 & \3 \\
\5 & \4 & \1 & \0 & \3 & \2 \\
\2 & \3 & \4 & \5 & \0 & \1 \\
\3 & \2 & \5 & \4 & \1 & \0 \\
\hline
\end{array}
\qquad
{7{:}}~
\begin{array}{|cc|ccccc|}
\hline
{\0}  & {\1}  & \2  & \3&\4&\5&\6\\
{\1}  & {\0}  & \3  & \4&\5&\6&\2\\ \hline
\2  & \4  & {\0}  & \1&\multicolumn{1}{|c}{\6}&\3&\5\\
\cline{3-3}\cline{5-5}\3  & \5  & \6  & \multicolumn{1}{|c}{\0}&{\1}&\multicolumn{1}{|c}{\2}&\4\\
\cline{4-4}\cline{6-6}\4  & \6  & \5  & \2&\multicolumn{1}{|c}{\0}&{\1}&\multicolumn{1}{|c|}{\3}\\
\cline{5-5}\cline{7-7}\5  & \2  & \4  & \6&\3&\multicolumn{1}{|c}{\0}&{\1}\\
\cline{3-3}\cline{6-6}\6  & \3  & \1  & \multicolumn{1}{|c}{\5}&\2&\4&\multicolumn{1}{|c|}{\0}\\
\hline
\end{array}
\end{equation}

Now, the statement follows from the obvious fact that a superposition
of $\Omega $-closed $2$-qua\-si\-groups is an $\Omega $-closed $n$-qua\-si\-group.
\end{Proof}

The next obvious lemma is a suitable tool for obtaining a large number of $n$-qua\-si\-groups,
most of which are irreducible.
\begin{Lemma}[(switching subquasigroups)]\label{l:switch}
Let $q:\Sigma^n\to\Sigma$ be an $\Omega$-closed $n$-qua\-si\-group with a subquasigroup
$g:\Omega ^n\to\Omega $, $g=q|_{\Omega ^n}$, $\Omega\subset\Sigma$.
And let $h:\Omega ^n\to\Omega$ be another $n$-qua\-si\-group of order $|\Omega |$.
Then
\begin{equation}\label{e:switch}
f(\bar x)\eqdf \cases{ h(\bar x) & if $\bar x\in\Omega ^n$ \cr q(\bar x) & if $\bar x\not\in\Omega ^n$}
\end{equation}
is an $n$-qua\-si\-group of order $|\Sigma|$.
\end{Lemma}

\section[Irreducible $n$-quasigroups]{{\hspace{-1em}.~\,Irreducible $n$-quasigroups}}\label{s:irreduc}

\begin{Lemma}\label{l:irsub}
A subquasigroup of a reducible $n$-qua\-si\-group is reducible.
\end{Lemma}
\begin{Proof}
Let $f:\Sigma^n\to\Sigma$ be a reducible $\Omega$-closed $n$-qua\-si\-group.
Without loss of generality we assume that
$$f(x,\bar y,\bar z)\equiv h(g(x,\bar y),\bar z)$$
for some $(n-m+1)$- and $m$-qua\-si\-groups $h$ and $g$ where $1<m<n$.
Take $\bar o\in \Omega^{m-1}$ and $\theta\in \Omega^{n-m}$.
Then the quasigroups $h_0$, $g_0$, and $\delta$ defined by (\ref{e:7})
are $\Omega$-closed. Therefore, the representation (\ref{e:8})
proves that $f|_{\Omega^n}$ is reducible.
\end{Proof}

\begin{Theorem}\label{th:irr}
For each $n\geq 3$ and $k\geq 4$ there exists an irreducible $n$-qua\-si\-group
of order~$k$.
\end{Theorem}
\begin{Proof}
(*) First we consider the case $n\geq 4$.
By Lemma~\ref{l:subq} we can construct
a reducible $n$-qua\-si\-group $q:\{0,\ldots,k-1\}^n\to\{0,\ldots,k-1\}$ of order $k$
with a subquasigroup $g:\{0,1\}^n\to\{0,1\}$ of order $2$.
Let $h:\{0,1\}^n\to\{0,1\}$ be the $n$-qua\-si\-group of order $2$ different from $g$;
and let $f$ be defined by (\ref{e:switch}).
By Theorem~\ref{th:vosst} with $\bar o=(2,\ldots,2)$,
the $n$-qua\-si\-group $f$ is irreducible.

(**) $n=3$, $k=4,5,6,7$.
In each of these cases we will construct an irreducible $3$-qua\-si\-group $f$,
omitting the verification, which can be done,
for example, using the formulas (\ref{e:7}), (\ref{e:8}).
Let quasigroups $q_{4,2}$, $q_{5,2}$, $q_{6,2}$, and $q_{7,2}$
be defined by the value arrays (\ref{e:57}).
For each case $k=4,5,6,7$ we define the ternary quasigroup
$q(x_1,x_2,x_3)\eqdf q_{k,2}(q_{k,2}(x_1,x_2),x_3)$,
which have the subquasigroup $q|_{\{0,1\}^3}(x_1,x_2,x_3)=x_1+x_2+x_3\bmod 2$.
Using (\ref{e:switch}), we replace this subquasigroup by
the ternary quasigroup $h(x_1,x_2,x_3)=x_1+x_2+x_3+1\bmod 2$.
The resulting ternary quasigroup $f$ is irreducible.

(***)  $n=3$, $8\leq k <\infty$.
Using Lemma~\ref{l:subq}, Lemma~\ref{l:switch}, and (**),
we can easily construct a ternary quasigroup of order $k\geq 8$ with
an irreducible subquasigroup of order $4$.
By Lemma~\ref{l:irsub}, such quasigroup is irreducible.

(****) The case of infinite order. Let $q:\Sigma_\infty^n\to\Sigma_\infty$ be an $n$-qua\-si\-group of infinite order $K$
and $g:\Sigma^n\to\Sigma$ be any irreducible $n$-qua\-si\-group of finite order (say, $4$).
Then, by Lemma~\ref{l:irsub}, their direct product
$g{\times} q:(\Sigma\times \Sigma_\infty)^n\to(\Sigma\times \Sigma_\infty)$
defined as
$$ g{\times} q\,\left([x_1,y_1],\ldots,[x_n,y_n]\vphantom{0^0}\right)
\eqdf
\left[g(x_1,\ldots,x_n),q(y_1,\ldots,y_n)\vphantom{0^0}\right]$$
is an irreducible $n$-qua\-si\-group of order $K$.
\end{Proof}

\begin{Note}
  Using the same arguments, it is easy to construct for any $n\geq 4$ and $k\geq 4$
   an irreducible $n$-qua\-si\-group of order $k$
  such that fixing one argument (say, the first) by (say) $0$ leads to
  an $(n-1)$-qua\-si\-group that is also irreducible.
  This simple observation naturally blends with the following context.
  Let $\kappa(q)$ be the maximal number such that
  there is an irreducible $\kappa(q)$-qua\-si\-group that can be obtained
  from $q$ or one of its inverses by fixing $n-\kappa(q)>0$ arguments.
  In this remark we observe that
  (for any $n$ and $k$ when the question is nontrivial)
  there is an irreducible
  $n$-qua\-si\-group $q$ with $\kappa(q)=n-1$.
  In \cite{Kro:n-2} for $k\vdots 4$ and even $n\geq 4$ an $n$-qua\-si\-group with
  $\kappa(q)=n-2$ is constructed.
  In \cite{Kro:n-3,KroPot:k2} it is shown that $\kappa(q)\leq n-3$
  (if $k$ is prime then $\kappa(q)\leq n-2$)
  implies that $q$ is reducible.
\end{Note}

\section[On the number of $n$-quasigroups, I]{{\hspace{-1em}.~\,On the number of $n$-quasigroups, I}}\label{s:number}
We first consider a simple bound on the number of $n$-qua\-si\-groups of composite order.
\begin{Prop}\label{p:composite}
The number $|Q(n,sr)|$ of $n$-qua\-si\-groups of composite order $sr$ satisfies
\begin{equation}\label{eq:composite}
|Q(n,sr)|\geq |Q(n,r)|\cdot|Q(n,s)|^{r^n}> |Q(n,s)|^{r^n}.
\end{equation}
\end{Prop}
\begin{Proof}
Let $g:Z_r^n \to Z_r$ be an arbitrary $n$-qua\-si\-group of order $r$;
and let $\omega\langle\cdot\rangle $ be an arbitrary function from $Z_r^n$ to the set $Q(n,s)$ of all
$n$-qua\-si\-groups of order $s$.
It is straightforward that the following function
is an $n$-qua\-si\-group of order $sr$:
$$
f(z\ij{1}{n}) \eqdf g\left(y\ij{1}{n}\vphantom{0^0}\right)\cdot s
+\omega \left\langle y\ij{1}{n}\vphantom{0^0}\right\rangle
( x\ij{1}{n}) \quad\mbox{ where } y_i \eqdf \lfloor z_i/s\rfloor, \quad x_i \eqdf z_i\bmod s
$$
$$
f(x_1,\ldots,x_n)=g\left(\lfloor x_1/s\rfloor,\ldots,\lfloor x_n/s\rfloor\vphantom{0^0}\right)\cdot s
+\omega \left\langle\lfloor x_1/s\rfloor,\ldots,\lfloor x_n/s\rfloor\vphantom{0^0}\right\rangle
( x_1\bmod s,\ldots, x_n\bmod s).
$$
Moreover, different choices of $\omega\langle\cdot\rangle$
result in different $n$-qua\-si\-groups.
So, this construction, which is known as the $\omega $-product of $g$,
obviously provides the bound (\ref{eq:composite}).
\end{Proof}
If the order is divided by $2$ or $3$
then the bound (\ref{eq:composite}) is the best known.
Substituting the known values $|Q(n,2)|=2$ and $|Q(n,3)|=3\cdot 2^n$, we get
\begin{Corol}\label{c:23}
  If $k\vdots 2$ then $|Q(n,k)|\geq 2^{(k/2)^n};$
  if $k\vdots 3$ then $|Q(n,k)|\geq (3\cdot 2^n)^{(k/3)^n}>2^{n(k/3)^n}$.
\end{Corol}
The next statement is weaker than the bound considered in the next section.
Nevertheless, it provides simplest arguments showing that the number of $n$-qua\-si\-group
of fixed order $k$ grows double-exponentially, even for prime $k\geq 8$.
The cases $k=5$ and $k=7$ will be covered in the next section.

\begin{Prop}\label{p:num1}
The number $|Q(n,k)|$ of $n$-qua\-si\-groups of order $k\geq 8$ satisfies
\begin{equation}\label{eq:bound0}
|Q(n,k)|\geq 2^{\lfloor k/4\rfloor^n}.
\end{equation}
\end{Prop}
\begin{Proof}
By Lemma~\ref{l:subq}, there is an $n$-qua\-si\-group of order $k$
with sub\-qua\-si\-group of order $2\lfloor k/4\rfloor$.
This sub\-qua\-si\-group can be switched (see Lemma~\ref{l:switch})
in $|Q(n,2\lfloor k/4\rfloor)|$ ways.
By Proposition~\ref{p:composite}, we have
$|Q(n,2\lfloor k/4\rfloor)|\geq |Q(n,2)|^{\lfloor k/4\rfloor^n}=2^{\lfloor k/4\rfloor^n}$.
Clearly, these calculations have sense only if $\lfloor k/4\rfloor>1$, i.\,e., $k\geq 8$.
\end{Proof}

\section[On the number of $n$-quasigroups, II]{{\hspace{-1em}.~\,On the number of $n$-quasigroups, II}}\label{s:comp}

In this section we continue using the same general switching principle as in previous ones:
independent changing the values of $n$-quasigroups on disjoint subsets of $\Sigma^n$.
We improve the lower bound in the cases when the order is not divided by $2$ or $3$;
in particular, we establish a double-exponential lower bound on the number of $n$-qua\-si\-groups
of orders $5$ and $7$.

We say that a nonempty set $\Theta\subset\Sigma^n$
is an \emph{$ab$-component} or a \emph{switching component} of an $n$-qua\-si\-group
$q$ iff
\begin{itemize}
  \item[(\textbf{a})] $q(\Theta)=\{a,b\}$ and
  \item[(\textbf{b})]
the function $q\Theta:\Sigma^n\to\Sigma$ defined as follows is an $n$-qua\-si\-group too:
  $$
  q\Theta(\bar x)\eqdf
  \cases{
   q(\bar x) & if $\bar x\not\in\Theta$\cr
   b & if $\bar x\in\Theta$ and $q(\bar x)=a$\cr
   a & if $\bar x\in\Theta$ and $q(\bar x)=b$.
  }
  $$
\end{itemize}

For example, $\{(0,0),(0,1),(1,0),(1,1)\}$ and $\{(2,2),(2,3),(3,3),(3,4),(4,2),(4,4)\}$
are $01$-components in (\ref{e:57}.$5$).

\begin{Note}
From some point of view, it is naturally to require also $\Theta$ to be inclusion-minimal, i.e.,
(\textbf{c}) $\Theta$ does not have a nonempty proper subset that satisfies (a) and (b).
Although in what follows all $ab$-components satisfy (c), formally we do not use it.
\end{Note}

\begin{Lemma}\label{l:countComp}
  Let an $n$-qua\-si\-group $q$ have $s$ pairwise disjoint switching components
  $\Theta_1$, \ldots, $\Theta_s$
  $($note that we do not require them to be $ab$-components for common $a$, $b)$.
  Then $|Q(n,|\Sigma|)|\geq 2^s$.
\end{Lemma}
\begin{Proof}
Indeed, denoting $q\Theta^0 \eqdf q$ and $q\Theta^1 \eqdf q\Theta$,
we have $2^s$ distinct $n$-quasigroups $q\Theta_1^{t_1}...\Theta_s^{t_s}$,
$(t_1,\ldots,t_s)\in\{0,1\}^s$.
\end{Proof}

\subsection[The order 5]{\hspace{-1em}.~\,The order 5}

In this section, we consider the $n$-quasigroups of order $5$, the only case,
when the other our bounds do not guarantee the double-exponential growth of the number
of $n$-quasigroups as $n\to\infty$. Of course, the way that we use for the order $5$
works for any other order $k>3$, but the bound obtained is
worse than (\ref{eq:composite}) provided $k$ is composite,
worse than (\ref{eq:bound0}) provided $k\geq 8$, and
worse than (\ref{eq:bound7}) provided $k\geq 6$.
The bound is based on the following straightforward fact:

\begin{Lemma}\label{l:for5}
Let $\{0,1\}^n$ be a $01$-component of an $n$-qua\-si\-group $q$.
For every $i\in\{1,\ldots,n\}$ let $q_i$ be an $n_i$-qua\-si\-group and
let $\Theta_i$ be its $01$-component. Then
$\Theta_1\times\ldots\times\Theta_n$ is a $01$-com\-po\-nent of the
$(n_1+\ldots +n_n)$-qua\-si\-group
$$
f(x_{1,1},...,x_{1,n_1},x_{2,1},\ldots,x_{n,n_n}) \eqdf q(q_1(x_{1,1},...,x_{1,n_1}),\ldots,q_n(x_{n,1},...,x_{n,n_n}))
.
$$
\end{Lemma}

For a quasigroup $q:\Sigma^2\to \Sigma$ denote $q^1\eqdf q$,
$q^2(x_1,x_2,x_3) \eqdf q(x_1,q^1(x_2,x_3))$, \ldots,
$q^i(x_1,x_2,\ldots,x_{i+1}) \eqdf q(x_1,q^{i-1}(x_2,\ldots,x_{i+1}))$.

\begin{Prop}\label{p:or5}
If $n=3m$ then $|Q(n,5)|\geq 2^{3^{m}};$
if $n=3m+1$ then $|Q(n,5)|\geq 2^{4\cdot 3^{m-1}};$
if $n=3m+2$ then $|Q(n,5)|\geq 2^{2\cdot 3^{m}}$.
Roughly, for any $n$ we have
$$|Q(n,5)|> 2^{3^{n/3-0.072}}>e^{e^{{\ln 3 \over 3}n - 0.44}}.$$
\end{Prop}
\begin{Proof}
  Let $q$ be the quasigroup of order $5$ with value table (\ref{e:57}.$5$).
Then

(*) $q$ has two disjoint $01$-components
$D_0\eqdf\{(0,0),(0,1),(1,0),(1,1)\}$ and
$D_1\eqdf\{(2,2),(2,3),(3,3),(3,4),(4,2),(4,4)\}$;

(**) $q^2$ has three mutually disjoint $01$-components
$T_0\eqdf 
\{0,1\}\times D_0$,
$T_1\eqdf 
\{0,1\} \times D_1$,
and
$T_2\eqdf \{(x_1,x_2,x_3)|q^2(x_1,x_2,x_3)\in\{0,1\}\}\setminus (T_0 \cup T_1)$;

(***) $\{0,1\}^{m+1}$ is a $01$-component of $q^m$.

By Lemma~\ref{l:for5},
\begin{itemize}
  \item[i.] the $3m$-quasigroup defined as the superposition
$$q^{m-1}\big(q^2(\cdot,\cdot,\cdot),\ldots,q^2(\cdot,\cdot,\cdot)\big)$$
has $3^m$ components $T_{t_1}\times\ldots\times T_{t_m}$, $({t_1},\ldots,{t_m})\in\{0,1,2\}^m$;
  \item[ii.]
the $3m+1$-quasigroup defined as the superposition
$$q^{m}\big(q^2(\cdot,\cdot,\cdot),\ldots,q^2(\cdot,\cdot,\cdot),q(\cdot,\cdot),q(\cdot,\cdot)\big)$$
has $3^{m-1}4$ components $T_{t_1}\times\ldots\times T_{t_{m-1}}\times D_{t_m}\times D_{t_{m+1}}$,
$({t_1},\ldots,{t_{m+1}})\in\{0,1,2\}^{m-1}\times\{0,1\}^2$;
  \item[iii.]
the $3m+2$-quasigroup defined as the superposition
$$q^{m}\big(q^2(\cdot,\cdot,\cdot),\ldots,q^2(\cdot,\cdot,\cdot),q(\cdot,\cdot)\big)$$
has $3^{m}2$ components $T_{t_1}\times\ldots\times T_{t_{m}}\times D_{t_{m+1}}$,
$({t_1},\ldots,{t_{m+1}})\in\{0,1,2\}^{m}\times\{0,1\}.$
\end{itemize}
By Lemma~\ref{l:countComp}, the theorem follows.
\end{Proof}
\begin{Note} If, in the proof, we consider the superposition
$q^{n/2}\big(q(\cdot,\cdot),\ldots,q^2(\cdot,\cdot)\big),$
then we obtain the bound $|Q(n,5)|\geq 2^{2^{n/2}}$ for even $n$,
which is worse because ${\ln 2 \over 2}<{\ln 3 \over 3}$.
\end{Note}

\subsection[The case of order $>6$]{\hspace{-1em}.~\,The case of order $\geq 7$}
In this section, we will prove the following:

\begin{Prop}\label{p:or7} The number $|Q(n,k)|$ of $n$-qua\-si\-groups
$\{0,1,\ldots,k-1\}^n\to \{0,1,\ldots,k-1\}$ satisfies
\begin{equation}\label{eq:bound7}
  |Q(n,k)|\geq 2^{\lfloor k/2\rfloor \lfloor k/3\rfloor ^{n-1}}
  > e^{e^{\ln \lfloor k/3\rfloor n + \ln \lfloor k/2\rfloor - \ln \lfloor k/3\rfloor -0.37}}
  > e^{e^{\ln \lfloor k/3\rfloor n +0.038}}
  .
\end{equation}

\end{Prop}
Note that this bound has no sense if $k<6$; and it is weaker than (\ref{eq:composite})
if $k\vdots 2$ or $k\vdots 3$. The proof is based on the following straightforward fact:
\begin{Lemma}\label{l:pow}
 Let $\{c,d\}\times\{e,f\}$ be an $ab$-com\-po\-nent of a quasigroup $g$. Then

 {\rm(a)} $\{a,b\}\times\{e,f\}$ is a $cd$-com\-po\-nent of the quasigroup ${g}^-$ defined by
 $g(x,y)=z \Leftrightarrow {g}^-(z,y)=x;$

  {\rm(b)} if $\{a_1,b_1\}\times\ldots\times\{a_n,b_n\}$ is an $ef$-com\-po\-nent of
  an $n$-qua\-si\-group $q$, then $\{c,d\}\times\{a_1,b_1\}\times\ldots\times\{a_n,b_n\}$
  is an $ab$-com\-po\-nent of
  the $(n+1)$-qua\-si\-group defined as the superposition $g(\cdot,q(\cdot,\ldots,\cdot))$.
\end{Lemma}

\begin{Proof}[of Proposition~\ref{p:or7}]
Taking into account Corollary~\ref{c:23},
it is enough to consider only the cases of odd $k \not\equiv 0 \bmod 3$.
Moreover, we can assume that $k>6$ (otherwise the statement is trivial).

Define the $2$-quasigroup $q$ as
\begin{eqnarray*}
q(2j,i)&\eqdf&i+3j\bmod k ;\\
q(2j+1,i)&\eqdf&\pi(i)+3j\bmod k;\\
q(2\lfloor k/3\rfloor+j,i)&\eqdf&\tau(i)+3j\bmod k;
\qquad j=0,\ldots,\lfloor k/3\rfloor-1,
\qquad i=0,\ldots,k-1
\end{eqnarray*}
where $\pi$, $\tau$, and the remaining values of $q$ are defined by the following
value table (the fourth row is used only for the case $k \equiv 2 \bmod 3$):
$$
\begin{array}{r||ccccccccccc}
i\hphantom{)}:       &  0 & 1 & \multicolumn{1}{|c}2 & \makebox[3.5ex][c]{3} & \multicolumn{1}{|c}{\makebox[3.5ex][c]4} & \ldots & \multicolumn{1}{|c}{k{-}5} & k{-}4 & \multicolumn{1}{|c}{k{-}3} & k{-}2 & k{-}1 \\
\pi(i):  &  1 & 0 & \multicolumn{1}{|c}3 & 2 & \multicolumn{1}{|c}{5} & \ldots & \multicolumn{1}{|c}{k{-}4} & k{-}5 & \multicolumn{1}{|c}{k{-}2} & k{-}1 & k{-}3 \cr
\hline
\tau(i): &k{-}1& \multicolumn{1}{|c}2 & 1 & \multicolumn{1}{|c}4 & 3 & \multicolumn{2}{|c|}{\ldots} & k{-}3 & \multicolumn{1}{c|}{k{-}4} &  0  & k{-}2 \cr
\hline
\hline
q(k{-}2,i):&k{-}3&k{-}2&k{-}1& 0 & 1 & \multicolumn{2}{c}{\ldots} & k{-}7 & k{-}6 & \multicolumn{1}{|c}{k{-}4} & k{-}5 \cr
q(k{-}1,i):&k{-}2&k{-}1& 0 & 1 & 2 & \multicolumn{2}{c}{\ldots} & k{-}6 & k{-}5 & \multicolumn{1}{|c}{k{-}3} & k{-}4 \cr
\end{array}
$$
In what follows, the tables illustrate the cases $k=7$ and $k=11$.
$$
\def\arraystretch{0.8}
\arraycolsep=0.6ex
\def\rsrs#1{\raisebox{-1.5pt}{{#1}}}
\def\g#1{\textcolor[gray]{0.6}{ #1}}
\def\0{\rsrs{0}}
\def\1{\rsrs{1}}
\def\2{\rsrs{2}}
\def\3{\rsrs{3}}
\def\4{\rsrs{4}}
\def\5{\rsrs{5}}
\def\6{\rsrs{6}}
\def\7{\rsrs{7}}
\def\8{\rsrs{8}}
\def\9{\rsrs{9}}
\def\a{\rsrs{\makebox[1ex][c]{1\hspace{-0.2ex}0}}}
k=7{:}~
\begin{array}{|cc|cc|ccc|}
\hline
\0&\1&\2&\3&\4&\5&\6\\
\1&\0&\3&\2&\5&\6&\4\\
\hline
\3&\4&\5&\6&\0&\1&\2\\
\4&\3&\6&\5&\1&\2&\0\\
\hline
\g\6&\multicolumn{1}{c}{\g\2}&\g\1&\multicolumn{1}{c}{\g\4}&\g\3&\g\0&\g\5\\
\g\2&\multicolumn{1}{c}{\g\5}&\g\4&\multicolumn{1}{c}{\g\0}&\g\6&\g\3&\g\1\\
\g\5&\multicolumn{1}{c}{\g\6}&\g\0&\multicolumn{1}{c}{\g\1}&\g\2&\g\4&\g\3\\
\hline
\end{array}
\qquad\qquad
k=11{:}~
\begin{array}{|cc|cc|cc|cc|ccc|}
\hline
\0&\1&\2&\3&\4&\5&\6&\7&\8&\9&\a\\
\1&\0&\3&\2&\5&\4&\7&\6&\9&\a&\8\\
\hline
\3&\4&\5&\6&\7&\8&\9&\a&\0&\1&\2\\
\4&\3&\6&\5&\8&\7&\a&\9&\1&\2&\0\\
\hline
\6&\7&\8&\9&\a&\0&\1&\2&\3&\4&\5\\
\7&\6&\9&\8&\0&\a&\2&\1&\4&\5&\3\\
\hline
\multicolumn{11}{|c|}{}
\end{array}
$$
For each
$j=0,\ldots,\lfloor k/3\rfloor-1$
and
$i=0,\ldots,\lfloor k/2\rfloor-2$
the set $\{2j,2j+1\}\times\{2i,2i+1\}$ is a $(2i+3j \bmod k)(2i+3j+1 \bmod k)$-com\-po\-nent
of such $q$. By Lemma~\ref{l:pow}(a), for the same pairs $i,j$
the set $\{2i+3j \bmod k,2i+3j+1 \bmod k\}\times\{2i,2i+1\}$ is a $(2j)(2j+1)$-com\-po\-nent of
$g\eqdf q^-$; moreover, we can observe that for each $j$ there is one more ``non-square''
$(2j)(2j+1)$-com\-po\-nent of $g$ which is disjoint with all considered ``square'' components,
see the following examples
(we omit the analytic description; indeed, we can ignore this component if we do not care about the constant in the bound
$e^{e^{\ln \lfloor k/3 \rfloor n + const}}$).
$$
\def\arraystretch{0.8}
\arraycolsep=0.6ex
\def\rsrs#1{\raisebox{-1.5pt}{{#1}}}
\def\g#1{\textcolor[gray]{0.7}{ #1}}
\def\0{\rsrs{0}}
\def\1{\rsrs{1}}
\def\2{\rsrs{2}}
\def\3{\rsrs{3}}
\def\4{\rsrs{4}}
\def\5{\rsrs{5}}
\def\6{\rsrs{6}}
\def\7{\rsrs{7}}
\def\8{\rsrs{8}}
\def\9{\rsrs{9}}
\def\a{\rsrs{\makebox[1ex][c]{1\hspace{-0.2ex}0}}}
k=7{:}~
\begin{array}{|cc|cc|ccc|}
\hline
\0&\1&\g\6&\g\5&\2&\g\4&\3\\
\1&\0&\g\4&\g\6&\3&\2&\g\5\\
\cline{1-4}
\g\5&\g\4&\0&\1&\g\6&\3&\2\\
\cline{1-2}\cline{5-7}
\2&\3&\1&\0&\g\4&\g\5&\g\6\\
\cline{3-7}
\3&\2&\g\5&\g\4&\0&\g\6&\1\\
\cline{1-4}
\g\6&\g\5&\2&\3&\1&\0&\g\4\\
\g\4&\g\6&\3&\2&\g\5&\1&\0\\
\hline
\end{array}
\qquad\qquad
\def\arraystretch{0.55}
\arraycolsep=0.4ex
\def\rsrs#1{\raisebox{-1.5pt}{{#1}}}
\def\0{\rsrs{$\scriptstyle 0$}}
\def\1{\rsrs{$\scriptstyle 1$}}
\def\2{\rsrs{$\scriptstyle 2$}}
\def\3{\rsrs{$\scriptstyle 3$}}
\def\4{\rsrs{$\scriptstyle 4$}}
\def\5{\rsrs{$\scriptstyle 5$}}
\def\6{\rsrs{$\scriptstyle 6$}}
\def\7{\rsrs{$\scriptstyle 7$}}
\def\8{\rsrs{$\scriptstyle 8$}}
\def\9{\rsrs{$\scriptstyle 9$}}
\def\a{\rsrs{\makebox[1ex][c]{$\scriptstyle 1\hspace{-0.2ex}0$}}}
k=11{:}~
\begin{array}{|cc|cc|cc|cc|ccc|}
\hline
  \0&\1&\g\a&\g\9&\5&\4&\g\8&\g\7&\2&\g\6&\3\\
\cline{5-8}
  \1&\0&\g\6&\g\a&\g\9&\g\8&\4&\5&\3&\2&\g\7\\
\cline{1-4}
\g\7&\g\6&\0&\1&\g\a&\g\9&\5&\4&\g\8&\3&\2\\
\cline{1-2}\cline{7-11}
  \2&\3&\1&\0&\g\6&\g\a&\g\9&\g\8&\4&\g\7&\5\\
\cline{3-6}
  \3&\2&\g\7&\g\6&\0&\1&\g\a&\g\9&\5&\4&\g\8\\
\cline{1-4}
\g\8&\g\7&\2&\3&\1&\0&\g\6&\g\a&\g\9&\5&\4\\
\cline{1-2}\cline{5-11}
  \4&\5&\3&\2&\g\7&\g\6&\0&\1&\g\a&\g\8&\g\9\\
\cline{3-6}
  \5&\4&\g\8&\g\7&\2&\3&\1&\0&\g\6&\g\9&\g\a\\
\cline{1-4}\cline{7-11}
\g\9&\g\8&\4&\5&\3&\2&\g\7&\g\6&\0&\g\a&\1\\
\cline{5-8}
\g\a&\g\9&\5&\4&\g\8&\g\7&\2&\3&\1&\0&\g\6\\
\cline{3-6}
\g\6&\g\a&\g\9&\g\8&\4&\5&\3&\2&\g\7&\1&\0\\
\hline
\multicolumn{1}{c}{}
\end{array}
$$
By induction, using Lemma~\ref{l:pow}(b), we derive that for every
$j_1,\ldots,j_{n-1}\in\{0,\ldots,\lfloor k/3\rfloor-1\}$ and $i\in\{0,\ldots,\lfloor k/2\rfloor-2\}$
the set
$$
\begin{array}{r@{}r@{~}r@{\,}c@{\,}l}
  \{&2j_2+3j_1 \bmod k,&2j_2+3j_1+1 \bmod k\}&\times&\\
  &\ldots\\
  \{&2j_{n-1}+3j_{n-2} \bmod k,&2j_{n-1}+3j_{n-2}+1 \bmod k\}&\times&\\
  \{&2i+3j_{n-1} \bmod k,&2i+3j_{n-1}+1 \bmod k\}&\times&\{2i,2i+1\}
\end{array}
$$
is a $(2j_1)(2j_1+1)$-com\-po\-nent of the $n$-qua\-si\-group $g^{n-1}$.
Also, for every such $j_1,\ldots,j_{n-1}$ there is one more
$(2j_1)(2j_1+1)$-com\-po\-nent of $g^{n-1}$, which is generated by the ``non-square''
$(2j_{n-1})(2j_{n-1}+1)$-com\-po\-nent of $g$.
In summary, $g^{n-1}$ has at least
$\lfloor k/3\rfloor^{n-1}\lfloor k/2\rfloor$ pairwise disjoint switching components.
By Lemma~\ref{l:countComp}, the theorem is proved.
\end{Proof}

Summarizing Corollary~\ref{c:23},
Propositions~\ref{p:or5} and~\ref{p:or7},
we get the following theorem.

\begin{Theorem}\label{th:Num}
  Let a finite set $\Sigma$ of size $k>3$ be fixed. The number $|Q(n,k)|$
  of $n$-quasigroups $\Sigma^n\to\Sigma$ satisfies the following$:$

{\rm (a)}  If $k$ is even, then $|Q(n,k)|\geq 2^{(k/2)^n}$.

{\rm (b)} If $k$ is divided by $3$, then $|Q(n,k)|\geq 2^{n(k/3)^n}$.

{\rm (c)} If $k=5$, then $|Q(n,k)|\geq 2^{3^{n/3-c}}$ where $c<0.072$ depends on $n\bmod 3$.

{\rm (d)} In all other cases, $|Q(n,k)|\geq 2^{1.5\lfloor k/3\rfloor ^n}$.
\end{Theorem}


\providecommand\href[2]{#2} \providecommand\url[1]{\href{#1}{#1}}

\end{document}